\documentclass[11pt]{amsart}
\usepackage[T1]{fontenc}
\usepackage{lmodern,microtype}
\usepackage[letterpaper,hmargin=1.16in,top=1.00in,bottom=1.03in]{geometry}
\usepackage{xcolor}
\definecolor{inkblue}{RGB}{29,63,93}
\usepackage[colorlinks=true,allcolors=inkblue,pdftitle={Area bounds for the genus and total curvature of simple pseudoholomorphic curves}]{hyperref}
\hypersetup{pdfauthor={Zhe Zhang}}
\newtheorem{theorem}{Theorem}[section]
\newtheorem{lemma}[theorem]{Lemma}
\newtheorem{proposition}[theorem]{Proposition}
\newtheorem{corollary}[theorem]{Corollary}
\theoremstyle{remark}
\newtheorem{example}[theorem]{Example}
\newtheorem*{remark}{Remark}
\newcommand{\II}{\mathrm{II}}
\newcommand{\PD}{\operatorname{PD}}
\newcommand{\Sec}{\operatorname{Sec}}
\newcommand{\im}{\operatorname{Im}}
\newcommand{\re}{\operatorname{Re}}
\newcommand{\reg}{\mathrm{reg}}
\newcommand{\darg}{\mathrm d\arg}
\numberwithin{equation}{section}
\title[Area bounds for simple pseudoholomorphic curves]{Area bounds for the genus and total curvature of simple pseudoholomorphic curves}
\author{Zhe Zhang}
\address{School of Arts and Sciences, North China Institute of Aerospace Engineering,
  Langfang, Hebei 065099, China}
\email{\href{mailto:zhangzhe@nciae.edu.cn}{zhangzhe@nciae.edu.cn}}
\date{November 27, 2026}
\begin{document}
\begin{abstract}
We prove a quadratic area bound for the genus of a closed simple pseudoholomorphic curve in a compact almost Hermitian manifold. The coefficients depend only on the ambient geometry, and embedded curves already attain quadratic growth. The same bound holds for the sum of the genus and the number of unordered self-intersection pairs with distinct complex tangent lines. The proof computes the self-intersection of the diagonal using an angular form on the product of the domain with itself. We also bound the total branch order and the integral of the squared second fundamental form. For a fixed symplectic form and background metric, these estimates are uniform over compatible almost complex structures with bounded $C^1$ norm. In real dimension six we recover the finiteness theorem of Doan and Walpuski~\cite[Theorem~1.6]{DW} by stable-map compactness. For primitive Calabi--Yau classes we obtain a quadratic upper bound in symplectic area for the genera supporting nonzero Gopakumar--Vafa (BPS) invariants, answering Question~1.14 of Doan and Walpuski~\cite{DWCount} in this case.
\end{abstract}
\maketitle

\section{Introduction}
Question~1.14 in Doan and Walpuski~\cite{DWCount} asks for a bound on the Castelnuovo number of a Calabi--Yau class $A$ in terms of $A$ and the geometry of the ambient closed symplectic six-manifold. The qualitative vanishing theorem \cite[Corollary~1.11]{DWCount} establishes that these invariants vanish at sufficiently high genus for primitive $A$. Corollary~\ref{cor:primitiveBPS} below gives a quadratic vanishing threshold in $\langle[\omega],A\rangle$ for all such classes. The underlying estimate applies to every closed simple pseudoholomorphic curve in a compact almost Hermitian manifold. It bounds the genus by a quadratic function of area, with a constant determined by the ambient geometry.

Let $u:(C,j)\to (X,J,g)$ be a nonconstant pseudoholomorphic map from a smooth closed connected Riemann surface. We fix a smooth compact almost Hermitian target and write $\eta(v,w)=g(Jv,w)$ for its fundamental form. The induced area form $dA=u^*\eta$ extends smoothly across the zeros of $du$ and has total mass $E$. On $C_{\reg}=\{du\ne0\}$ we use the induced metric, write $\II$ for the second fundamental form and take $H=\operatorname{tr}\II$. At a self-intersection these quantities are evaluated on each sheet in the domain.

For an immersion, the Gauss equation and the Gauss--Bonnet theorem \cite[Chapters~8--9]{Lee} give
\begin{equation}\label{eq:gaussintro}
 \int_C |\II|^2\,dA=\int_C\bigl(2\Sec_X(du(TC))+|H|^2\bigr)\,dA-4\pi\chi(C).
\end{equation}
An upper bound on the ambient sectional curvature and on $|H|$ controls the positive part of the intrinsic curvature. The term $-|\II|^2/2$ in the Gauss equation leaves open the lower bound for $\int K_C\,dA$ needed to control the genus. For an embedded curve with $\langle c_1(TX,J),u_*[C]\rangle=0$, the complex normal bundle expresses the same problem as bounding $\langle c_1(N_C),[C]\rangle=2g(C)-2$.

A curve is \emph{simple} if it does not factor through a holomorphic branched cover of degree greater than one. At a critical point $p$ let $b_p$ be the vanishing order of $du$. The local expansion in \cite[Appendix~B.2.3, equation~(B.15)]{Wendl} gives $u(z)-u(p)=az^{b_p+1}+O(|z|^{b_p+2})$ with $a\ne0$ in smooth target coordinates complex linear at $u(p)$. Put $B(u)=\sum_p b_p$. Take $b_p=0$ at a regular point. Denote by $N_{\rm tr}(u)$ the number of unordered pairs $\{p,q\}$ with $p\ne q$, $u(p)=u(q)$ and distinct complex tangent lines, using the limiting tangent line at a critical point. This definition also applies in higher codimension. If several sheets pass through the same point, we count each unordered pair with distinct tangent lines. The genus throughout is that of the smooth domain.

\begin{theorem}\label{thm:main}
There exists $K>0$ depending only on $(X,J,g)$ such that every closed connected simple $J$-holomorphic curve satisfies
\begin{equation}\label{eq:genus}
 g(C)+N_{\rm tr}(u)\leq\left(\frac{KE}{2\pi}-1\right)^2 .
\end{equation}
Let $S_+=\max\{0,\sup\Sec_X\}$, and let $H_*$ be a uniform bound for $|H|$ on the regular locus of such curves. Then
\begin{equation}\label{eq:curv}
 \int_{C_{\reg}}|\II|^2\,dA+4\pi B(u)+8\pi N_{\rm tr}(u)
 \leq (2S_++H_*^2)E+\frac{2K^2}{\pi}E^2-8KE .
\end{equation}
The same $K$ applies to embedded curves and to curves with singularities. Each distinct-line pair may in fact be counted with weight $(b_p+1)(b_q+1)$ in both estimates.
\end{theorem}

\noindent\textit{Sharpness of the quadratic growth rate.}
Let $Y$ be a projective K3 surface, $L$ an ample line bundle on $Y$, and $T$ an elliptic curve. Give $Y\times T$ a fixed product K\"ahler metric and write $\omega_Y$ for the K\"ahler form on $Y$. For large $n$, a general member $C_n$ of $|L^{\otimes n}|$ is smooth and connected \cite[II.8.18 and III.7.9.1]{Hartshorne}. Regard $C_n$ as an embedded curve in $Y\times\{t\}$ for a fixed $t\in T$. Its area and the adjunction formula \cite[II.8.20]{Hartshorne} give, with intersections taken on $Y$,
\[
 E_n=n[\omega_Y]\cdot c_1(L),\qquad 2g(C_n)-2=n^2 c_1(L)^2 .
\]
Since $C_n$ is minimal and the ambient sectional curvature is bounded, \eqref{eq:gaussintro} gives $\int_{C_n}|\II|^2\,dA=4\pi n^2 c_1(L)^2+O(n)$. Thus both the genus and the total squared second fundamental form grow quadratically in $E_n$.

\begin{remark}[Simplicity]
Double covers $C_g\to\mathbb P^1$ of a fixed projective line have constant area and arbitrarily large genus, with total branch order $2g+2$ by Riemann--Hurwitz \cite[IV, \S2]{Hartshorne}. The graph of the deck involution gives a second curve of coincidences in $C_g\times C_g$, so division by the diagonal leaves a positive-dimensional zero set.
\end{remark}

\begin{corollary}\label{cor:symp}
If $J$ is compatible with a symplectic form $\omega$ and $g=\omega(\cdot,J\cdot)$, one can take $H_*=0$ and $E=\langle[\omega],u_*[C]\rangle$. Thus both estimates apply in every represented class of a fixed almost K\"ahler manifold.
\end{corollary}

\begin{theorem}\label{thm:family}
Fix a compact symplectic manifold $(X,\omega)$, a smooth background metric $h$, and a number $M_1$. There exist $K,C_0>0$ depending only on $(X,\omega,h,M_1)$ with the following property. For every smooth $\omega$-compatible $J$ with $\|J\|_{C^1(h)}\le M_1$, every closed connected simple $J$-holomorphic curve satisfies \eqref{eq:genus} with this $K$ and
\begin{equation}\label{eq:familystatement}
 \frac12\int_{C_{\reg}} |\II_{g_J}|^2\,dA+4\pi B(u)+8\pi N_{\rm tr}(u)
 \le C_0 E+\frac{2K^2}{\pi}E^2-8KE .
\end{equation}
Here $g_J=\omega(\cdot,J\cdot)$, $dA=u^*\omega$, and $\II_{g_J}$ is the second fundamental form computed with the Levi--Civita connection of $g_J$. An area bound thus controls the genus, the number of distinct-line pairs, the total branch order, and the total squared second fundamental form.
\end{theorem}

The normal bundle of the diagonal $\Delta\subset C\times C$ is $TC$, so the self-intersection formula \cite[\S11]{BottTu} gives $\Delta^2=\langle e(TC),[C]\rangle=\chi(C)$. To estimate this self-intersection number, we pull the angular form of the ambient displacement back to $C\times C$. For an embedding the only pole is along $\Delta$. For a simple curve, the leading normal derivative can vanish on $\Delta$, and distinct source points can map to the same target point. The regularization in Section~\ref{sec:product} removes these additional singularities while retaining the phase of a local defining function for $\Delta$. The boundary of the diagonal tube gives the same Euler characteristic term in both cases. At a self-intersection with distinct tangent lines, the local contribution is the degree of a map between three-spheres times $4\pi^2$ for each ordered pair of preimages. For families, Section~\ref{sec:family} constructs a unitary connection with first-order control of its curvature and compares its projected tangent connection with the one induced by $g_J$.

The earlier finiteness theorem of Doan and Walpuski~\cite[Theorem~1.6]{DW} gives a qualitative genus bound for fixed $k$-rigid structures in dimension six under its Chern-number and divisibility hypotheses \cite[Remark~1.7]{DW}. Their proof uses compactness for pseudoholomorphic cycles \cite[Proposition~1.9]{DW}. Section~\ref{sec:finite} instead obtains the genus bound from Theorem~\ref{thm:main} and applies stable-map compactness twice. For the BPS application we use the primitive-class identification of BPS numbers with signed embedded-curve counts \cite[footnote~11]{Zinger}, as derived in \cite[Appendix~C]{DWCount}.

\section{Angular curvature and radial regularization}\label{sec:angular}
Equip $X\times X$ with $G=g\oplus g$ and $\mathcal J=J\oplus J$. Choose a unitary connection on $TX$ and pull it back by the first projection. Denote this connection by $\nabla$. In a fixed metric tube $U$ about the ambient diagonal, the section
\[
 D(x,y)=\exp_x^{-1}(y)\in T_xX
\]
is smooth, and $|D|=d_g(x,y)$. Hermitian pairings are conjugate-linear in the first entry. For any section $s$ of a Hermitian bundle with unitary connection define, where $s\ne0$,
\[
 \theta(s)=\frac{\im\langle s,\nabla s\rangle}{|s|^2}.
\]
In particular $\theta(hs)=\darg h+\theta(s)$ for a nonvanishing complex function $h$, where $\darg h=\im(h^{-1}dh)$.

We use a radial profile $f:[0,\infty)\to[0,1]$ that is smooth, nondecreasing, concave, linear near zero with $f(0)=0$, and identically one on $[1,\infty)$. Such a function is the normalized integral from $0$ to $t$ of a nonnegative, nonincreasing smooth function that is constant near zero and vanishes for $t\ge1$. Concavity gives
\begin{equation}\label{eq:profile}
 0\leq t f'(t)\leq f(t)\leq1 .
\end{equation}
Set $\theta_{\epsilon}(s)=f(|s|^2/\epsilon^2)\theta(s)$. As a form pulled back by a smooth section it extends smoothly across $s=0$, since near zero it is a constant times $\epsilon^{-2}\im\langle s,\nabla s\rangle$.

\begin{lemma}\label{lem:radial}
Let $(Y,I)$ be almost complex and let $s$ be a section of a Hermitian bundle on $Y$ with unitary connection of curvature $R$. For $v\in TY$ define
\[
 B_\pm(v)=\tfrac12\big(\nabla_v s\mp i\nabla_{Iv}s\big).
\]
At $s\ne0$ let $Q$ be the Hermitian projection onto $\mathbb C s$, set $P=1-Q$, and put $t=|s|^2/\epsilon^2$. Then
\begin{align}
 d\theta_\epsilon(s)(v,Iv)
 &=\frac{2}{|s|^2}\big[
   f(t)(|PB_+|^2-|PB_-|^2)
     +t f'(t)(|QB_+|^2-|QB_-|^2)\big] \notag\\
 &\hspace{1em}+f(t)\frac{\im\langle s,R(v,Iv)s\rangle}{|s|^2}.\label{eq:weighted}
\end{align}
The unregularized formula for $d\theta(s)$ is obtained by setting $f=1$ and $f'=0$. In either case the right-hand side is bounded below by
\begin{equation}\label{eq:uniformradial}
 -2\frac{|B_-(v)|^2}{|s|^2}-\lVert R(v,Iv)\rVert_{\rm op}.
\end{equation}
\end{lemma}
\begin{proof}
Write $V=\nabla_v s$, $W=\nabla_{Iv}s$. Differentiating the angular form of the unit section $s/|s|$ gives $(2\im\langle PV,PW\rangle+\im\langle s,R(v,Iv)s\rangle)/|s|^2$. The term parallel to the unit section drops out because its covariant derivative in that line is purely imaginary. The product rule adds $f'(t)dt\wedge\theta(s)$. Since $dt(v)=2\re\langle s,V\rangle/\epsilon^2$, this term is $2t f'(t)\im\langle QV,QW\rangle/|s|^2$. For either projection $T$, write $V=B_++B_-$ and $W=i(B_+-B_-)$. Thus $\im\langle TV,TW\rangle=|TB_+|^2-|TB_-|^2$. Equation~\eqref{eq:weighted} follows. The two projections are orthogonal, so \eqref{eq:profile} bounds the sum of the two negative terms by $2|B_-|^2/|s|^2$ in absolute value.
\end{proof}

Choose a radial cutoff $\rho(|D|)$ equal to one in an inner tube and zero near the outer edge of $U$. Extend $a_X=\rho(|D|)\theta(D)$ by zero outside $U$. On the ambient diagonal,
$\nabla_{(v_x,v_y)}D=v_y-v_x$, which is complex linear. Smoothness and compactness give a constant $M$ on the closed inner tube such that
\begin{equation}\label{eq:ambient}
 |(\nabla D)^{0,1}(v)|\le M|D||v|_G .
\end{equation}
Here $(\nabla D)^{0,1}(v)$ is $B_-(v)$ for $s=D$. Let $R_*$ bound $\|R(v,\mathcal Jv)\|_{\rm op}$ for $|v|_G=1$ there, and let $L$ bound the norm of $da_X$ on the cutoff annulus. Fix once and for all
\begin{equation}\label{eq:K}
 K_0=\max\{2M^2+R_*,L\},\qquad K>K_0 .
\end{equation}
By Lemma~\ref{lem:radial}, $da_X(v,\mathcal Jv)\ge-K_0|v|_G^2$ everywhere off the ambient diagonal. This is the constant obtained from the unregularized embedded-curve estimate.

\section{The intrinsic diagonal of a simple curve}\label{sec:product}
Write $Y=C\times C$, $I=j\oplus j$ and $F=u\times u$. Let $\Delta\subset Y$ and $\Delta_X\subset X\times X$ denote the respective diagonals. With $p_1,p_2$ the projections, put
\begin{equation}\label{eq:omega}
 \Omega=p_1^*dA+p_2^*dA,\qquad
 \Omega(v,Iv)=|dF(v)|_G^2,\qquad
 \int_Y\Omega^2=2E^2.
\end{equation}
This closed form can be degenerate over the critical points of $u$. In the inverse image of the ambient tube let $D_u=F^*D$ denote the displacement section in the pullback of $TX$. Away from $F^{-1}(\Delta_X)$ the ambient angular form pulls back to $F^*a_X$.

\subsection{Division by the diagonal}
Consider the exceptional set of ordered pairs
\begin{equation}\label{eq:exception}
 Z=\{(p,p):du_p=0\}\ \cup\ \{(p,q):p\ne q,\ u(p)=u(q)\}.
\end{equation}
Wendl~\cite[Theorem~A.1]{Wendl} gives finiteness of both the critical points and the double points of a closed simple curve. This makes $Z$ finite, including tangential intersections and points with several preimages. The local branch structure underlying this result is described by Micallef and White~\cite[Theorems~6.1--6.2]{MW}.

Choose disjoint coordinate neighborhoods of the points of $Z$, with closures in the inverse image of the inner tube. Near $(p,p)$ take $h=z_2-z_1$ in a common holomorphic coordinate on $C$. Near an off-diagonal pair take the neighborhood disjoint from $\Delta$ and put $h=1$. These are local representations of the canonical defining section of $\mathcal O(\Delta)$. Set $s=D_u/h$ where $h\ne0$. Holomorphicity of $h$ and complex linearity of $dF$ give
\begin{equation}\label{eq:division}
 (\nabla s)^{0,1}(v)=h^{-1}(\nabla D_u)^{0,1}(v),\qquad
 |(\nabla s)^{0,1}(v)|\le M|s|\,|dF(v)|_G .
\end{equation}
The factor $|h|^{-1}$ cancels against $|D_u|=|h||s|$. This cancellation retains the ambient estimate even in directions where $dF$ vanishes.

To see how the quotient behaves on the diagonal, work in a smooth unitary frame over the first variable and write $z_1=z$, $z_2=z+r e^{i\vartheta}$. Taylor expansion with integral remainder gives
\begin{equation}\label{eq:polar}
 D_u(z,z+r e^{i\vartheta})=r e^{i\vartheta}a(z)+r^2 b(z,r,\vartheta),
 \qquad a(z)=du_z(\partial_x).
\end{equation}
Here $b$ is smooth down to $r=0$ in the displayed polar variables. The leading term follows from complex linearity of $du$ and the identity derivative of the exponential map at zero. Thus $s$ extends continuously on the product, and smoothly in polar variables, with boundary value $a(z)$ independent of $\vartheta$. Complex linearity of the first jet alone does not imply smooth division in the Cartesian variables. For example, the smooth scalar function $h+\bar h^2$ has complex-linear derivative at zero, but its quotient by $h$ is not differentiable there. Only the polar extension and its tangential derivatives enter the proof. In particular the only zeros left, including the extended values on $\Delta$, are those in \eqref{eq:exception}.

\subsection{Removal of the residual zeros}
On a closed collar of the boundary of each chosen neighborhood, $|s|$ has a positive minimum after the continuous extension. Shrinking the neighborhoods if necessary, choose $\epsilon>0$ below all these minima. If $Z$ is empty no modification is needed. On the punctured product $Y\setminus\Delta$ define a one-form $\alpha$ by retaining $F^*a_X$ outside the neighborhoods and putting inside them
\begin{equation}\label{eq:alpha}
 \alpha=\darg h+\theta_\epsilon(s).
\end{equation}
The formulas agree on the collars, since there $f=1$, $\rho=1$ and $\theta(D_u)=\darg h+\theta(s)$. At an off-diagonal zero $h=1$ and $\theta_\epsilon(s)$ is smooth. Thus $\alpha$ is smooth on all of $Y\setminus\Delta$. Both the neighborhoods and $\epsilon$ may depend on the individual map.

\begin{proposition}\label{prop:lower}
The form $\beta=d\alpha$ satisfies $\beta(v,Iv)\ge-K_0\Omega(v,Iv)$ on $Y\setminus\Delta$, for the $K_0$ of \eqref{eq:K}.
\end{proposition}
\begin{proof}
Outside the neighborhoods the assertion is the pullback of the ambient estimate. Inside, $d(\darg h)=0$. Apply \eqref{eq:uniformradial} to $s=D_u/h$, use \eqref{eq:division}, and bound the pulled-back curvature by $R_*|dF(v)|_G^2$. This gives the same $2M^2+R_*$ as in the inner tube. At a regularized zero the inequality follows by continuity. The bound is uniform in $\epsilon$.
\end{proof}

\begin{example}\label{ex:axes}
For the flat cusp $u(z)=(z^2,z^3)$ one obtains
$s(z,w)=(z+w,z^2+zw+w^2)$ after division by $w-z$. Its only zero is $(0,0)$, and its restriction to the diagonal is $(2z,3z^2)$. Formula~\eqref{eq:alpha} regularizes the phase of this vector while retaining $\darg(w-z)$. 
For two coordinate axes meeting at distinct source points, $h=1$ and $D_u=(-z,w)$ on the product of the two parameter disks. The unregularized phase has a nonzero boundary term. On the unit sphere it restricts to the standard Hopf contact form $\lambda=\im(\bar z\,dz+\bar w\,dw)$ in the normalization of \cite[\S2]{AGZ}. With the complex orientation,
\[
 \int_{S^3}\theta(D_u)\wedge d\theta(D_u)
 =\int_{B^4}(d\lambda)^2=8\operatorname{vol}(B^4)=4\pi^2 .
\]
Here $\lambda$ on the right is the same smooth one-form on the ball. Dilation gives the same flux on arbitrarily small spheres. Formula~\eqref{eq:alpha} fills the ball smoothly and retains its boundary phase. Simply omitting the boundary of an isolated coincidence would lose this contribution.
\end{example}

\subsection{The boundary circle and its curvature}
Fix a smooth $j$-Hermitian metric on $C$, including at its critical points. Parametrize a tube about $\Delta$ by $(p,\xi)\mapsto(p,\exp_p\xi)$ for this auxiliary metric. Its radius-$r$ boundary is identified with the unit circle bundle $b:S(TC)\to C$ by $\xi=rv$. Write $Y_r$ for the complement of the open tube. In the order base, positive fiber circle, the boundary orientation of $Y_r$ is the negative of the circle-bundle orientation.

The restriction of $\alpha$ to the circle bundle converges with all tangential derivatives as $r\to 0$ to a circle connection $a_0$. Indeed, away from the modified charts the section $D_u/r$ tends smoothly to the nonzero vector $du_p(v)$, which is equivariant under rotation of $v$. Its angular form is invariant and takes the value one on the infinitesimal positive rotation. In a modified chart meeting $\Delta$, the function $h/r$ tends to the nonzero complex-linear evaluation $dz_p(v)$. Thus $\darg h$ on the boundary tends to $\darg(dz_p(v))$ with the same rotation and normalization properties. The quotient $s=(D_u/r)/(h/r)$ tends smoothly to $du_p(\partial_x)$, independent of $v$, as in \eqref{eq:polar}. Consequently the tangential limit of $\theta_\epsilon(s)$ there is the pullback of a smooth base one-form, including at $du_p=0$. The limits agree on collars. Also, the boundary restriction of $\Omega$ tends to $2b^*dA$.

The normalized angular form $a_0/(2\pi)$ has fiber integral one, and its derivative is the negative pullback of an Euler form \cite[\S\S6, 11]{BottTu}. With $c_1=[iR/(2\pi)]$, write $a_0=d\vartheta+A_0$ in a local unitary frame of $TC$. Then $da_0=b^*F_0$, where $F_0=dA_0$ and the line connection has curvature $iF_0$. Hence
\begin{equation}\label{eq:residue}
 \int_{S(T_pC)}a_0=2\pi,\qquad \int_C F_0=-2\pi\chi(C).
\end{equation}
In particular the regularization retains the diagonal class. More precisely, for every smooth closed two-form $\psi$ on $Y$, Stokes' theorem gives
\begin{equation}\label{eq:pairing}
 \lim_{r\to 0}\int_{Y_r}\beta\wedge\psi
 =-\int_{S(TC)}a_0\wedge b^*(\psi|_\Delta)
 =-2\pi\int_\Delta\psi .
\end{equation}
This identifies the pairing with $-2\pi\PD[\Delta]$, including its sign.

\subsection{The square and the genus}
Proposition~\ref{prop:lower} makes $\beta+K\Omega$ nonnegative on complex lines. Its square is a nonnegative four-form. Indeed, adding any positive multiple of the fundamental form of an auxiliary Hermitian metric gives a taming form. Each taming form induces the complex orientation by the straight-line homotopy to the auxiliary form. Letting the positive multiple decrease to zero gives the assertion also where $\Omega$ is degenerate.

On $Y_r$ the square equals
$d\big(\alpha\wedge(\beta+2K\Omega)\big)+K^2\Omega^2$.
Stokes' theorem and \eqref{eq:residue} now yield
\begin{align}
 0&\le\lim_{r\to 0}\int_{Y_r}(\beta+K\Omega)^2 \notag\\
  &=-\int_{S(TC)}a_0\wedge b^*(F_0+4K\,dA)+2K^2E^2 \notag\\
  &=4\pi^2\chi(C)-8\pi KE+2K^2E^2. \label{eq:stokes}
\end{align}
The normal line of $\Delta$ is $TC$ via $(v,w)\mapsto w-v$. Consequently $\Delta^2=\chi(C)$ and the last expression is the square of $K[\Omega]-2\pi\PD[\Delta]$ on $Y$. Its nonnegativity follows from the pointwise estimate and the boundary calculation. For embeddings, substituting $\chi(C)=2-2g(C)$ gives \eqref{eq:genus} immediately.

\subsection{The local contribution of distinct tangent lines}\label{sub:hopf}
We compute how much of the nonnegative square is concentrated near an ordered coincidence $a=(p,q)$ with $p\ne q$. Write $k=b_p+1$, $\ell=b_q+1$, taking $b_p=0$ at a regular point. The local expansion \cite[Appendix~B.2.3, equation~(B.15)]{Wendl} identifies the two limiting tangent lines with the spans of the leading Taylor coefficients, in target coordinates complex linear at $u(p)=u(q)$.

\begin{lemma}\label{lem:hopf}
If the two tangent lines at $a$ are distinct, there are neighborhoods $U_r$ shrinking to $a$ for which
\begin{equation}\label{eq:flux}
 \lim_{r\to 0}\int_{\partial U_r}F^*a_X\wedge d(F^*a_X)=4\pi^2 k\ell,
 \qquad
 \lim_{r\to 0}\int_{\partial U_r}F^*a_X\wedge\Omega=0 .
\end{equation}
Both integrals use the boundary orientation induced by the complex orientation of $U_r$.
\end{lemma}
\begin{proof}
Choose holomorphic coordinates at $p,q$ and dilate the unit ball in $\mathbb C^2$ by $\delta_r(z,w)=(r^{1/k}z,r^{1/\ell}w)$, defining $U_r=\delta_r(B^4)$. In a smooth unitary frame centered at the image point, Taylor expansion gives
\[
 r^{-1} D_u(\delta_r(z,w))\longrightarrow P(z,w)=-v z^k+v' w^\ell
 \quad\text{in }C^2(S^3).
\]
The leading terms of the two branches are $v z^k$ and $v'w^\ell$. The differential of the displacement takes their difference, and all higher Taylor terms carry a positive power of $r$ after division by $r$. The rescaled connection one-form tends to zero in $C^1$ on the sphere, since it is pulled back via the first map, which tends in $C^2$ to a constant. The complex independence of $v,v'$ gives $\min_{S^3}|P|>0$. It follows that the angular form converges in $C^1$ to the pullback of the standard unit angular form $\lambda$ by $P/|P|$. The dilated pullback of $\Omega$ tends to zero because both component maps converge in $C^1$ to a constant.

Identify $\operatorname{span}_{\mathbb C}(v,v')$ unitarily with $\mathbb C^2$. The degree of $P/|P|:S^3\to S^3$ is the local degree of $P$ at the origin. In the basis $(-v,v')$, a sufficiently small regular value with both coordinates nonzero has $k\ell$ preimages, all with positive Jacobian. The degree is therefore $k\ell$, regardless of the angle between the two lines. The flat computation in Example~\ref{ex:axes} gives $\int_{S^3}\lambda\wedge d\lambda=4\pi^2$, proving \eqref{eq:flux}.
\end{proof}

Let $\mathcal D$ be the set of unordered pairs of distinct preimages with distinct limiting complex tangent lines, and define
\begin{equation}\label{eq:weight}
 W(u)=\sum_{\{p,q\}\in\mathcal D}(b_p+1)(b_q+1).
\end{equation}
Each unordered pair gives two points in $C\times C$, even if other pairs have the same image in $X$. First choose $r$ so that the neighborhoods of all these points are disjoint and avoid $\Delta$. With $r$ fixed, choose the regularization scales in \eqref{eq:alpha} below the minima of $|D_u|$ on the boundaries. Then $\alpha=F^*a_X$ there. For a fixed smooth regularization the integral would tend to zero over a shrinking ball. The total integral in \eqref{eq:stokes} is independent of the regularization, and its integrand is nonnegative. On each ball, Stokes' theorem gives
\[
 \int_{U_r}(\beta+K\Omega)^2
 =\int_{\partial U_r}F^*a_X\wedge\bigl(d(F^*a_X)+2K\Omega\bigr)
       +K^2\int_{U_r}\Omega^2 .
\]
The boundary term involving $\Omega$ tends to zero by Lemma~\ref{lem:hopf}. The volume term tends to zero because $\Omega$ is smooth. Summing over the ordered pairs and letting $r\to 0$ gives
\begin{equation}\label{eq:budget}
 4\pi^2\chi(C)-8\pi KE+2K^2E^2\ \ge\ 8\pi^2 W(u).
\end{equation}
Since $\chi(C)=2-2g(C)$, this is equivalent to
\[
 g(C)+W(u)\le\left(\frac{KE}{2\pi}-1\right)^2.
\]
Now $W(u)\ge N_{\rm tr}(u)$ gives \eqref{eq:genus}. Each unordered pair of regular sheets with distinct tangent lines contributes one to $W(u)$. Pairs with the same tangent line are omitted from \eqref{eq:weight}.

\section{The branch contribution to total curvature}\label{sec:curvature}
In this section $K_C$ denotes the Gaussian curvature of $u^*g$ on $C_{\reg}$. The smooth mean-curvature estimate goes back to Gray~\cite[Proposition~5.5]{Gray}. Its formulation for pseudoholomorphic cycles is given in \cite[Theorem~2.36 in the preprint]{DIW}. The first identity below is the conical Gauss--Bonnet formula \cite[\S2, Proposition~1]{Troyanov} with divisor $\sum_p b_p\,p$. We include the boundary calculation to fix the sign and justify the curvature integrals for the induced metric. Both identities hold for nonconstant multiply covered maps.
\begin{proposition}\label{prop:branch}
For a nonconstant closed connected pseudoholomorphic map one has
\begin{align}
 \int_{C_{\reg}}K_C\,dA &= 2\pi\bigl(\chi(C)+B(u)\bigr),\label{eq:gb}\\
 \int_{C_{\reg}} |\II|^2\,dA+4\pi B(u)
 &=\int_{C_{\reg}}\bigl(2\Sec_X(du(TC))+|H|^2\bigr)\,dA
       -4\pi\chi(C).\label{eq:branchedgauss}
\end{align}
Both integrals on the left are finite, and $K_C$ is absolutely integrable.
\end{proposition}
\begin{proof}
First, $|H|\le H_*:=\|d\eta\|_{\rm comass,\infty}$ on $C_{\reg}$. On each immersed patch the area equals the integral of $\eta$, which has comass one. Under a compactly supported two-sided variation the difference has a minimum at zero. Its first derivative therefore vanishes, so every compactly supported normal variation field $V$ satisfies
\[
 -\int\langle H,V\rangle\,dA=\int u^*(\iota_V d\eta).
\]
Localizing this identity gives the pointwise bound.

Near a critical point of order $b$, choose a holomorphic source coordinate $z$ and smooth target coordinates with $u(0)=0$, $J(0)=i$ and $g(0)$ unitary. The first nonzero Taylor term is $az^{b+1}$ for $a\ne0$ by the expansion preceding equation~(B.15) in \cite[Appendix~B.2.3]{Wendl}. Taylor's theorem for the smooth map and its derivatives thus gives
\[
 u^*g=e^{2\lambda}|dz|^2,\qquad
 \lambda=b\log r+q(r,\vartheta),\qquad
 q=\log((b+1)|a|)+O(r),\qquad \partial_rq=O(1).
\]
For the last bound, $r^{-b}u_x$ is smooth in polar coordinates down to $r=0$ with nonzero leading vector $(b+1)ae^{ib\vartheta}$. On the positively oriented circle of radius $r$, the geodesic curvature satisfies
\[
 \int_{|z|=r} k_g\,ds
 =\int_0^{2\pi}(1+r\partial_r\lambda)\,d\vartheta
 =2\pi(b+1)+O(r).
\]
Each hole in the excised surface carries the opposite boundary orientation. If there are $N$ critical points, its Euler characteristic is $\chi(C)-N$. Applying Gauss--Bonnet there and shrinking the holes gives \eqref{eq:gb} as an improper integral.

On the regular locus the Gauss equation reads $2K_C=2\Sec_X(du(TC))+|H|^2-|\II|^2$. Its positive curvature part has finite integral because the area and the ambient terms are bounded. Together with the finite limit in \eqref{eq:gb} this proves absolute integrability. Integrating the Gauss equation gives \eqref{eq:branchedgauss} and finiteness of the second fundamental form integral.
\end{proof}

For a simple map insert \eqref{eq:budget} in \eqref{eq:branchedgauss} and bound the ambient terms by $(2S_++H_*^2)E$. This proves \eqref{eq:curv}. One can retain $8\pi W(u)$ on the left. Dropping $-8KE$ recovers the bound $(2S_++H_*^2)E+2K^2E^2/\pi$ for the entire left-hand side. For compatible symplectic data $d\eta=d\omega=0$ and $dA=u^*\omega$, which gives Corollary~\ref{cor:symp}. A double point with two immersed sheets changes neither the Euler characteristic of the smooth domain nor the local Gauss--Bonnet integrals. Only zeros of $du$ occur in $B(u)$.

\section{Uniform estimates for compatible structures with a first derivative bound}\label{sec:family}
Fix a compact symplectic manifold $(X,\omega)$, a smooth background metric $h$, and $M_1<\infty$. All ambient constants in this section depend only on these data. The hypothesis $\|J\|_{C^1(h)}\le M_1$ includes a bound on $J$ itself. Compatibility makes $g_J=\omega(\cdot,J\cdot)$ positive definite. Both $g_J$ and $g_J^{-1}$ are uniformly bounded relative to $h$, since the metric is the composition of the fixed symplectic isomorphism with $J$ and $J^{-1}=-J$. The first background derivative of $g_J$ is bounded as well.

\subsection{A unitary connection with uniformly bounded curvature}
Choose a fixed symplectic connection $D^0$ on $TX$, for example by patching the flat connections of Darboux charts with a partition of unity, and extend it complex linearly. The indefinite Hermitian form $\mathcal H(\xi,\zeta)=i\omega_{\mathbb C}(\bar\xi,\zeta)$ on $TX\otimes\mathbb C$ is parallel. For $P=(1-iJ)/2$ its range $V=T^{1,0}_J X$ is $\mathcal H$-positive and $\ker P=V^{\perp_{\mathcal H}}$. The map $\Phi_J(v)=(v-iJv)/\sqrt 2$ is a complex-linear isometry from $(TX,J,g_J)$ onto $V$ with its induced Hermitian metric. Project $D^0$ onto $V$ and transport the resulting unitary connection back to $TX$. Denote it by $D^J$. The Gauss formula for a projected connection \cite[\S45, equation~(45.8)]{Derdzinski} has the following expression in terms of $P$. Using $P(D^0P)P=0$ gives
\begin{align}
 D^J&=D^0-\tfrac12J(D^0J),\notag\\
 R^{P D^0}|_V&=\big[P R^0 P+P(D^0P\wedge D^0P)P\big]|_V .\label{eq:projection}
\end{align}
The second formula follows by applying $P D^0$ twice to a section $P s=s$ and replacing $(1-P)D^0s$ by $(D^0P)s$. In particular the norms of $R^{D^J}$ and the torsion $T^{D^J}$ are uniformly bounded. Also the difference $D^J-\nabla^{g_J}$ from the Levi--Civita connection is bounded since $g_J$ is uniformly elliptic with bounded first derivatives. These estimates require only $D^0J$.

\subsection{Uniform choice of the diagonal constant}
Use the displacement $D_h(x,y)=(\exp^h_x)^{-1}y$ and a cutoff in the fixed distance $d_h(x,y)$. Measure $D_h$ in $g_J$ and use the connection $D^J$ in its angular form. Here the fixed tube and cutoff are independent of $J$. To obtain \eqref{eq:ambient} uniformly, let $\tau_{yx}$ be $h$-parallel transport from $y$ to $x$ and set $r=d_h(x,y)$. In the fixed tube,
\[
 (D^0 D_h)(v_x,v_y)=\tau_{yx}v_y-v_x+O(r)(|v_x|_h+|v_y|_h).
\]
The fixed smooth coefficients have this first-order bound since they agree on the diagonal. Replacing $D^0$ by $D^J$ adds $O(r)|v_x|_h$. Furthermore $|\tau_{yx}J_y-J_x\tau_{yx}|_h\le c r$ by the first derivative bound. Applying the last two estimates to $(v_x,v_y)$ and $(J_xv_x,J_yv_y)$ gives the required $O(r)|v|$ bound on the anti-complex-linear part. Norm equivalence identifies $r$ with $|D_h|_{g_J}$ up to uniform factors. The annular cutoff estimate follows directly from the angular derivative formula in Lemma~\ref{lem:radial}: $|D_h|$ there has a uniform positive lower bound, and $|D^J D_h|$, $|R^{D^J}|$ and the derivative of the cutoff are uniformly bounded. Thus one $K$ works in \eqref{eq:budget} throughout the family, with $dA=u^*\omega$ and $E=\int_Cu^*\omega$.

\subsection{The curvature of the tangent line}
We estimate the Riemannian $\II_{g_J}$ itself. On the immersed locus let $L=du(TC)$ with the connection obtained by projecting $D^J$ to this complex line. Write its curvature as $iF_L$. The normal projection of $D^J$ on tangent fields defines the bilinear form $\mathcal A(X,Y)=(D^J_XY)^\perp$; here tangent vectors are identified with their images by $du$. In a local oriented unit frame $e,Je$ put
\[
 A=\mathcal A(e,e),\quad V=\mathcal A(Je,e),\quad
 q=\im\langle e,R^{D^J}(e,Je)e\rangle .
\]
Since $D^J$ preserves $J$, $\mathcal A(X,Je)=J\mathcal A(X,e)$. The projected torsion identity gives $JA-V=T^{D^J}(e,Je)^\perp$. Differentiating the real connection form $\im\langle e,D^Je\rangle$ as in Lemma~\ref{lem:radial} gives
$F_L(e,Je)=q+2g_J(JA,V)$. Consequently
\begin{equation}\label{eq:linecurv}
 |\mathcal A|^2
 =2(|A|^2+|V|^2)
 =2F_L(e,Je)-2q+2|T^{D^J}(e,Je)^\perp|^2 .
\end{equation}

To integrate it, compare the two projected connections on $L$. On $C_{\reg}$ the projection of $\nabla^{g_J}$ is the intrinsic Levi--Civita connection; it is unitary because an oriented Riemannian surface preserves rotation by $\pi/2$ under parallel transport. The two line connections differ by $i\gamma$ with $|\gamma|_{u^*g_J}\le c$. Their real curvatures therefore satisfy $F_L=-K_C\,dA+d\gamma$. At a critical point of order $b$, the induced length of the radius-$r$ circle is $O(r^{b+1})$ by the local expansion in Proposition~\ref{prop:branch}. Thus the $\gamma$-integrals on all excision boundaries tend to zero. Branched Gauss--Bonnet now gives
\begin{equation}\label{eq:lineint}
 \int_{C_{\reg}}F_L=-2\pi\big(\chi(C)+B(u)\big).
\end{equation}
These boundary limits are taken for each smooth map separately. The lower bound on $F_L(e,Je)$ from \eqref{eq:linecurv} and the finite integral in \eqref{eq:lineint} also give absolute integrability.

For a uniform $c_0$, integration of \eqref{eq:linecurv} yields $\int|\mathcal A|^2\,dA\le -4\pi(\chi(C)+B(u))+c_0 E$. Pointwise, $|\II_{g_J}-\mathcal A|\le c_1$ by the connection difference estimate. Applying $\tfrac12|\II_{g_J}|^2\le |\mathcal A|^2+c_1^2$ and inserting \eqref{eq:budget} gives, with $C_0=c_0+c_1^2$,
\begin{equation}\label{eq:familystrong}
 \frac12\int_{C_{\reg}}|\II_{g_J}|^2\,dA+4\pi B(u)+8\pi W(u)
 \le C_0 E+\frac{2K^2}{\pi}E^2-8KE .
\end{equation}
This proves Theorem~\ref{thm:family}. Although the integral of $K_C$ was evaluated for each individual metric, no bound for its pointwise size or for $\Sec(g_J)$ enters the constants.

\begin{example}
On the fixed symplectic torus with angular coordinates $x,y$ and $\omega=dx\wedge dy$, the metrics $g_n=f_n^{-1}dx^2+f_n\,dy^2$ with $f_n(x)=1+(2n)^{-1}\sin(nx)$ determine compatible structures with bounded $C^1$ norms. Their Gaussian curvatures are $-f_n''/2$ and have supremum $n/4$. This also gives examples in any higher even dimension by taking a product with a fixed flat symplectic torus.
\end{example}

\section{Finiteness by stable-map compactness}\label{sec:finite}
The proof uses the normal rescaling of \cite[\S3]{DW} and Gromov compactness \cite[\S1]{Gromov} in the stable-map formulation of Ivashkovich and Shevchishin~\cite[Theorem~1.1]{IS}. Theorem~\ref{thm:main} supplies a fixed-genus subsequence for both applications of compactness, first in $X$ and then in the normal bundle of an embedded limiting curve.

\begin{theorem}[Doan--Walpuski {\cite[Theorem~1.6]{DW}}]\label{thm:finite}
Let $(X,J,g)$ be a compact almost Hermitian six-manifold and let $J$ be $k$-rigid in the sense of \cite[Definition~2.10]{DW}, for $k\in\mathbb N\cup\{\infty\}$. Suppose $A\in H_2(X;\mathbb Z)$ satisfies $\langle c_1(TX,J),A\rangle=0$ and has divisibility at most $k$. For every $\Lambda>0$ there are finitely many simple closed $J$-holomorphic curves in class $A$ of area at most $\Lambda$, up to reparametrization.
\end{theorem}

\begin{proof}
Assume that $u_\nu:C_\nu\to X$ form an infinite sequence of inequivalent curves as in the statement. The index formula \cite[equation~(2.5)]{DW} gives unmarked index $2\langle c_1(TX,J),A\rangle=0$ in real dimension six. The embedding condition in \cite[Definition~2.10(3)]{DW} therefore applies. The right-hand side of \eqref{eq:genus} is bounded on $[0,\Lambda]$, so a subsequence has fixed genus. The compactness theorem \cite[Theorem~1.1]{IS} gives a connected stable limit of total class $A$, with Hausdorff convergence of the images. The limit has a nonconstant component. Otherwise the images lie eventually in a single coordinate ball carrying an exact symplectic form taming $J$, whose integral over each map would be both positive and zero.

Using the factorization of closed pseudoholomorphic maps into a branched cover and a simple map \cite[Appendix~A.1]{Wendl}, write $B_i$ for the distinct simple images and $m_i$ for their total covering degrees. Additivity gives
\[
 0=2\langle c_1(TX,J),A\rangle
   =\sum_i m_i\,2\langle c_1(TX,J),[B_i]\rangle .
\]
Parts~(2) and~(3) of the rigidity definition \cite[Definition~2.10]{DW} make each summand nonnegative and ensure that distinct simple index-zero images are embedded and disjoint. Constant components attach at a single image point, so connectedness gives one image $B$ and $A=m[B]$. Divisibility gives $m\le k$ for finite $k$. For $k=\infty$ the finite sum of covering degrees already supplies a finite $m$.

For large $\nu$ the curves lie in a fixed tubular neighborhood of $B$. Write $v_\nu:C_\nu\to N_B$ for the maps in tubular coordinates and $\pi:N_B\to B$ for the projection. Then $\pi v_\nu$ has degree $m$. Indeed, after collapsing the shrinking circles to nodes as in \cite[Definition~1.6]{IS}, the projected limit is uniformly close to $\pi v_\nu$ and hence homotopic to it by short geodesics in $B$. The collapse sends $[C_\nu]$ to the sum of the component classes, whose projected degrees sum to $m$. This computes the degree in the tube itself.

Denote fiber dilation by $\delta_\lambda$. Discard the single possible curve with image $B$. With $|v_\nu|$ the fiber norm, set
\[
 \lambda_\nu=2\max_{C_\nu}|v_\nu|\longrightarrow0,\qquad
 \widetilde u_\nu=\delta_{\lambda_\nu}^{-1}v_\nu .
\]
The rescaled curves lie in the closed disk bundle of radius $1/2$ and meet its boundary. By \cite[Lemma~3.10]{DW}, the rescaled structures $J_\nu$ converge smoothly on compact sets to the structure $J_B$ of the normal Cauchy--Riemann operator. The projection is $J_B$-holomorphic, and each fiber has its usual complex structure \cite[Lemma~3.6]{DW}.

To control rescaled energies choose an area form $\sigma$ on $B$, a unitary connection on $N_B$, and the fiberwise Liouville form $\vartheta$ extended by zero on the horizontal distribution. On the closed radius-one disk bundle
\[
 \omega_N=\pi^*\sigma+t\,d\vartheta
\]
tames $J_B$ for all sufficiently small $t>0$. Indeed, $d\vartheta$ is the standard positive fiber form on the vertical block, has zero mixed block, and has bounded horizontal block. In this splitting $J_B$ has diagonal blocks $j_B,i$ and a bounded horizontal-to-vertical block. For a tangent vector with horizontal and vertical lengths $h,v$, the mixed error is at most $t c h v$, and the horizontal error is at most $t c h^2$. Young's inequality absorbs these into the positive base term and $t$ times the positive fiber term. Fixing one such $t$, smooth convergence gives a uniform taming bound for $J_\nu$. Since dilation preserves the projection, Stokes' theorem gives
\[
 \int_{C_\nu}\widetilde u_\nu^*\omega_N=m\int_B\sigma.
\]
This bounds the energies in a fixed metric on the disk bundle.

Apply the varying-target-structure case of \cite[Theorem~1.1]{IS} on a neighborhood of the closed radius-$1/2$ disk bundle. The genus is fixed and the images lie in a common compact set. The limit meets the sphere bundle of radius $1/2$. A nonconstant component meets it there as well, since a connected union of constant components and its attaching nodes have the same image. Each nonconstant component projects to $B$ with positive degree, since a holomorphic map from a closed curve to a single vector-space fiber is constant. These degrees sum to $m$. Factor a component meeting the radius-$1/2$ boundary through its simple image map $w$. Then $\varphi=\pi w$ has $1\le\deg\varphi\le m$, and $w$ gives a nonzero section of $\varphi^*N_B$ in the kernel of the pulled-back normal operator \cite[Proposition~3.12]{DW}. This contradicts $k$-rigidity of $B$.
\end{proof}

\subsection{The BPS Castelnuovo bound for primitive classes}
\begin{corollary}[Primitive case of {\cite[Question~1.14]{DWCount}}]\label{cor:primitiveBPS}
On a closed symplectic six-manifold $(X,\omega)$, fix a smooth metric $h$, a smooth $\omega$-compatible structure $J_*$, and $M_1>\|J_*\|_{C^1(h)}$. Let $K$ be the constant of Theorem~\ref{thm:family} for these data. For every primitive class $A\in H_2(X;\mathbb Z)$ with $\langle c_1(TX),A\rangle=0$ and $E_A=\langle[\omega],A\rangle>0$, one has
\begin{equation}\label{eq:primitiveBPS}
 \mathrm{BPS}_{A,g}(X,\omega)=0
 \qquad\text{for every integer }g>\left(\frac{K E_A}{2\pi}-1\right)^2 .
\end{equation}
\end{corollary}
\begin{proof}
The strict $C^1$ bound defines a nonempty open set of smooth compatible structures containing $J_*$. The residuality statement \cite[Proposition~2.37(1)]{DWCount} gives a $J$ in this set for which the primitive signed embedded counts $n_{A,g}(X,\omega)$ are defined and independent of the choice of regular $J$ by \cite[Theorem~1.5(1)]{DWCount}. For primitive Calabi--Yau classes the comparison with Gromov--Witten invariants \cite[footnote~11]{Zinger} identifies $n_{A,g}(X,\omega)$ with $\mathrm{BPS}_{A,g}(X,\omega)$. Appendix~C of \cite{DWCount} derives this identity using the comparison theorem \cite[Theorem~1.2]{Zinger} and the local contribution formula \cite[\S2.3]{Pandharipande}. Each counted curve has area $E_A$. Theorem~\ref{thm:family} makes the counted moduli space empty for every $g$ in \eqref{eq:primitiveBPS}. The bound on $J$ fixes $K$ before $A$ is chosen.
\end{proof}

\end{document}